\documentclass{article}
\usepackage{amssymb,amsmath,amsthm}
\usepackage{amsfonts}
\usepackage{amscd}
\newtheorem{lem}{Lemma}
\newtheorem{thm}{Theorem}
\begin{document}

\title{About one singularity of the property of disturbed
differential systems, that equivalents to linear differential
equations}
\author{I.~Kopshaev
\thanks{Institute of mathematics of NAS of KAZAKHSTAN, 125 Pushkina str.,
050010 Almaty, KAZAKHSTAN email: {\em{kopshaev@math.kz}}}}

\maketitle
\date

\begin{abstract}
This paper is devoted to the study of the family of morphisms of
vector bundle, defined by the systems equivalence to linear
differential equations. It is proved, that the specified families
of the morphisms of vector bundle is not saturated.
\end{abstract}

\section{Introduction}
In what follows, I consider the families of morphisms of vector
bundle \cite{Mill1}, that defines by the linear systems of
ordinary differential equations.

In \cite{Mill2} proved, that the families of the morphisms of
vector bundles, defines by the arbitrary systems of differential
equations are saturated. The property of saturation, for the first
time, was defined in \cite{Mill3}. It turned out \cite{Rahim1},
\cite{Rahim2}, that using of this property allows to define the
conditions of stability of Liapunov exponent of morphisms of
vector bundles.

The aim of this paper is the study of the property of saturation
conformably to the families of morphisms, that defines by the
systems equivalence to linear differential equations.

\section{Preliminaries}
Let $(E,p,B)$ be a vector bundle with the fibre $R^{n}$ and the
base $B$ (where $B$ is a full metric space). Fix on $(E,p,B)$ some
Riemannian metric (\cite{Husemoller}, P. 58-59).

Consider the homomorphism of group $Z$ (group $R$) into group of
isomorphisms of vector bundle $(E,p,B)$. Remind, that this means
following \cite{Mill1}: for all $t \in Z$ ($t \in R$) is given
$X^{t}$ - homomorphism $E$ into $E$ and $\chi^{t}$ - homomorphism
$B$ into $B$ such that, $pX^{t}=\chi^{t}p$; for all $b \in B$ the
contraction $X^{t}[b]$ of map $X^{t}$ on fibre $p^{-1}(b)$ is
linear map $p^{-1}(b) \to p^{-1}(\chi^{t}b)$; for all $t,s \in Z$
(correspondingly $R$) have a places equalities $X^{t+s}=X^{t}
\cdot X^{s}$, $\chi^{t+s}=\chi^{t} \cdot \chi^{s}$.

The image of point $t$ for this homomorphism denotes as $(X^{t},
\chi^{t})$, instead of $X^{1}$ writes $X$, instead of $\chi^{1}$ -
$\chi$.

Suppose, that there exists the function $a(\cdot): B \to R^{+}$,
that satisfies the equality $a(\chi^t b)=a(b)$ for all $b \in B$
and for all $t \in Z$ (correspondingly $t \in R$) such that, for
all $t \in N$ (correspondingly $t \in R^{+}$) has place inequality
\begin{equation}
\label{kop_eq1} \max(\| X^{t}[b] \|, \|X^{-t}[b] \|) \leq e^{t
\cdot a(b)}
\end{equation}
(the norm of a linear map of fibre into a fibre defines in the
ordinary way by using norms in the fibres, induced as mentioned
above fixed Riemannian metric of vector bundle $(E,p,B)$).

Assign for all $m \in N$
$$
(X(m), \chi(m)) \stackrel{\rm def}{=} (X^{m}, \chi^{m}).
$$

Obtained in this way the family of morphisms $(X^{m}, \chi^{m})$
$m \in N$ of vector bundle $(E,p,B)$, that satisfies the following
condition: there exists the function $a(\cdot):B \to R^{+}$ such
that, for all $b \in B$, $m \in N$ has place inequality
\begin{equation}
\label{kop_eq2} \max(\| X(m,b) \|, \| \left[ X(m,b) \right]^{-1}
\|) \leq e^{m \cdot a(b)},
\end{equation}
where as $X(m,b)$ denotes the contraction of map $X(m)$ on fibre
$p^{-1}(b)$; thus, if $X(m)=X^{m}$, then $X(m,b)=X^{m}[b]$.
Moreover, such defined family of morphisms $(X(m), \chi(m))$ ($m
\in N$) satisfies the conditions a) - c):

a) $(X, \chi)$ - isomorphism of vector bundle $(E,p,B)$;

b) for all $m \in N$ have places equalities
$$
X(m) = X^{m}, \; \chi(m) = \chi^{m};
$$
c) there exists the function $a(\cdot):B \to R^{+}$ such that,
$a(\chi^{m}b)=a(b)$ for all $b \in B$, $m \in Z$, and so that, for
all $b \in B$ has place inequality
\begin{equation}
\label{kop_eq3} \max(\| X[b] \|, \| \left[ X[b] \right]^{-1} \|)
\leq e^{a(b)}.
\end{equation}
The verification of all this claims is trivial.

Remind the definition of saturated families of morphisms, given in
\cite{Mill3}:

{\bf Definition}\,{\it The family of morphisms
$(X(m),\chi(m)):(E,p,B) \to (E,p,B)$ ($m \in N$), that satisfies
the conditions a) - c), is said to be saturated, if for any point
$b \in B$ such that $\chi^{m}b \neq b$ for $\forall m \neq 0$, for
$\forall \varepsilon > 0$, for any basis $\{ \xi_{1}, \ldots,
\xi_{n} \}$ of vector space $p^{-1}(b)$ and for all neighbourhoods
$U(\xi_{i})$ of points $\xi_{i}$ ($i \in \{1,\ldots,n\}$) (in the
space $E$) there exists $\delta> 0$ such that for $\forall
\overline{t} \in N$ and for all nonsingular linear operators
$$
Y_{m}: p^{-1}(\chi^{m-1}b) \to p^{-1}(\chi^{m}b)
$$
($m \in \{1,\ldots,\overline{t}\}$), that satisfy for $\forall m
\in \{1,\ldots,\overline{t}\}$ the inequality
\begin{equation}
\label{kop_eq4} \| Y_{m}(X[\chi^{m-1}b])^{-1} - E \| +
\|X[\chi^{m-1}b] Y^{-1}_{m} - E\| < \delta
\end{equation}
there exists the point $b' \in B$ such that
\begin{equation}
\label{kop_eq5} d_{B}(b',b) < \varepsilon,
\end{equation}
and for $\forall m \in \{0,\ldots,\overline{t}\}$ there exist
isomorphisms of fibres (as Cartesian spaces)
$$
\psi_{m}: p^{-1}(\chi^{m}b') \to p^{-1}(\chi^{m}b),
$$
in which connection fulfiled following conditions:

i) $ \psi_{0}^{-1} \xi_{i} \in U(\xi_{i})$ for $\forall i \in \{
1, \ldots, n\}$;

ii) for $\forall m \in \{1,\ldots,\overline{t}\}$ diagram
\begin{equation}
\label{kop_eq6}
\begin{CD}
p^{-1}(\chi^{m-1}b') @>{X[\chi^{m-1}b']}>>  p^{-1}(\chi^{m}b')  \\
@VV{\psi_{m-1}}V                            @VV{\psi_{m}}V   \\
p^{-1}(\chi^{m-1}b) @>{Y_{m}}>>            p^{-1}(\chi^{m}b)
\end{CD}
\end{equation}
is commutative.}

Consider the family of morphisms defined in \cite{Rahim3}.

Let the space of all linear systems of differential equations
\begin{equation}
\label{kop_eq7} \dot{x}=A(t) \cdot x, \; x \in R^{n},
\end{equation}
such that $A(\cdot):R \to \mathrm{Hom}(R^{n},R^{n})$ - continuous
map, for which $\sup\limits_{t \in R}\| A(t)\|< + \infty$, endow
with a structure of the metric space by define distance by formula
$$
d(A_{1},A_{2}) = \sup\|A_{1}(t) - A_{2}(t)\|.
$$
(here point $\dot{x}=A_{i}(t) \cdot x$ of space denotes by
$A_{i}$). Obtained metric space $M_{n}$ is full.

Let
\begin{equation}
\label{kop_eq8} B \stackrel{\rm def}{=} M_{n}, \, E \stackrel{\rm
def}{=} B \times R^{n}, \, p \stackrel{\rm def}{=} pr_{1},
\end{equation}
where $pr_{1}$ - projection of production $B \times R^{n}$ on the
first efficient.

Thus, the trivial vector bundle $(E,p,B)$ is defined.

Let for all $t \in R$
\begin{equation}
\label{kop_eq9} X^{t}(A,x)= (\chi(t)A, \mathfrak{X}(t,0,A)x),
\end{equation}
\begin{equation}
\label{kop_eq10} \chi^{t}A(\cdot)= A(t+(\cdot)),
\end{equation}
where $A \in B$, $x \in R^{n}$, $\mathfrak{X}(\Theta,\tau,A)$ -
the Cauchy matrix of a system $\dot{x}=A(t) \cdot x$.

For all $t \in R$ have places formulas:
$$
X^{t} \cdot X^{-t}=X^{0}=1_{E}, \; X^{-t} \cdot X^{t}=X^{0}=1_{E},
$$
$$
\chi^{t} \cdot \chi^{-t}=\chi^{0}=1_{B}, \; \chi^{-t} \cdot
\chi^{t}=\chi^{0}=1_{B},
$$
i.e. $(X^{-t},\chi^{-t})$ - morphism of vector bundle $(E,p,B)$,
inverse for $(X^{t},\chi^{t})$, therefor, $(X^{t},\chi^{t})$ -
isomorphism of vector bundle $(E,p,B)$.

Thus, the homomorphism of group $R$ into group of isomorphisms of
vector bundle $(E,p,B)$ is defined. The image of point $t \in R$
for this homomorphism is $(X^{t},\chi^{t})$.

The function $a(\cdot):B \to R^{+}$ defines by formula
$$
a(A) \stackrel{\rm def}{=} \sup\limits_{t \in R} \|A(t)\|.
$$
We have
$$
a(\chi^{t}A) = a(A)= \sup\limits_{t \in R} \|A(t)\|.
$$
The correctness for all $b \in B$, $t \in R^{+}$ of inequality
(\ref{kop_eq1}) follows from (\ref{kop_eq9}), in view of
well-known inequality
$$
\| \mathfrak{X}(t,0,A) \| \leq e^{|\int_{0}^{t} \| A(s) ds\| |}
\leq e^{|t \cdot a(A)|},
$$
to which Cauchy matrix of system  $\dot{x}=A(t) \cdot x$ satisfies
for all $A \in B$, $t \in R$.

Let
$$
(X(m), \chi(m)) \stackrel{\rm def}{=} (X^{m},\chi^{m}).
$$
We obtain the family of morphisms
$$
(X(m), \chi(m)): (E,p,B) \to (E,p,B),
$$
($m \in N$), that satisfies the conditions a) - c).

In \cite{Mill2} proved, that the family of morphisms
(\ref{kop_eq8}) is saturated.

\section{Auxiliary results}
Consider the differential equation
\begin{equation}
\label{kop_eq11} \ddot{y}=a(t) \cdot y, \; y \in R,
\end{equation}
where $a(\cdot): R \to R$ - continues map with $\sup\limits_{t \in
R}|a(t)| < \infty$.

The equation (\ref{kop_eq11}) equivalents to the system of
differential equations
\begin{equation}
\label{kop_eq12} \dot{x}=A(t)x, \; x \in R^{2}, \; A(t)=
\begin{pmatrix}
0          & 1 \\
a(t)       & 0
\end{pmatrix}.
\end{equation}
Consider also the differential equations
\begin{equation}
\label{kop_eq13} \ddot{y}=[a(t)+b_{\varepsilon}(t)] \cdot y, \; y
\in R,
\end{equation}
where $b_{\varepsilon}(\cdot): R \to R$ - continues map with
$\sup\limits_{t \in R}|b_{\varepsilon}(t)| < \varepsilon$
$(\varepsilon > 0)$.

The equation (\ref{kop_eq13}) equivalents to system of
differential equations
\begin{equation}
\label{kop_eq14} \dot{x}=[A(t)+B_{\varepsilon}(t)] \cdot x, \; x
\in R^{2},
\end{equation}
where matrix $B_{\varepsilon}(t)$ has representation and
estimation
\begin{equation}
\label{kop_eq15} B_{\varepsilon}(t)=\left(
\begin{array}{cc}
0          & 0 \\
b_{\varepsilon}(t)       & 0
\end{array}\right), \; \sup\limits_{t \in R}\|B_{\varepsilon}(t)\| <
\varepsilon.
\end{equation}

\begin{lem}
\label{kop_lm1} If the differential equation (\ref{kop_eq11}) is
given and $\mathfrak{X}(\Theta, \tau,A)$ - the Cauchy matrix of a
equivalent system (\ref{kop_eq12}), then for $\forall \varepsilon
> 0$ there exists $\delta > 0$ such that for $\forall \overline{t}
\in N$ there always exists nonsingular linear operators $W_{m}:R^2
\to R^2$ ($m \in \{1, \ldots, \overline{t}\}$) that satisfy for
$\forall m \in \{1, \ldots, \overline{t}\}$ the inequality
\begin{equation}
\label{kop_eq16} \| W_{m}[\mathfrak{X}(m,m-1,A)]^{-1} - E \| +
\|\mathfrak{X}(m,m-1,A) W^{-1}_{m} - E\| < \delta
\end{equation}
for which there no exists continues map $$
A_{\varepsilon}(\cdot):[0,\overline{t}] \to \mathrm{Hom}(R^2,R^2)
$$ that satisfies the condition:

1)$\sup\limits_{t \in [0,\overline{t}]}\|A_{\varepsilon}(t)-A(t)\|
< \varepsilon$;

2) $\mathfrak{X}(m,m-1,A_{\varepsilon})=W_{m}$ for $\forall m \in
\{1, \ldots, \overline{t}\}$,

\noindent where $\mathfrak{X}(m,m-1,A_{\varepsilon})$ - Cauchy
matrix of system (\ref{kop_eq14}), that equivalents to equation
(\ref{kop_eq13}).
\end{lem}

{\it Proof.} We argue by contradiction. Then for $\forall
\varepsilon>0$ there exists $\delta
> 0$ such that for $\forall \overline{t} \in N$ and for all
nonsingular linear operators $W_{m}:R^2 \to R^2$ ($m \in \{1,
\ldots, \overline{t}\}$), that satisfy the inequality
(\ref{kop_eq16}), there exists continues map
$A_{\varepsilon}(\cdot):[0,\overline{t}] \to
\mathrm{Hom}(R^2,R^2)$, that satisfies the conditions 1), 2).

Let $\mathfrak{X}(t,s,A_{\varepsilon})$ be Cauchy matrix of system
(\ref{kop_eq14}). Then for $\forall t,s \in R$ has place equality:
$$
\dot{\mathfrak{X}}(t,s,A_{\varepsilon}) = A_{\varepsilon}(t) \cdot
\mathfrak{X}(t,s,A_{\varepsilon})=[A(t)+B_{\varepsilon}(t)] \cdot
\mathfrak{X}(t,s,A_{\varepsilon}).
$$
If we multiply last equality from the right by
$\mathfrak{X}^{-1}(t,s,A_{\varepsilon})$ , then obtain equality
\begin{equation}
\label{kop_eq17} B_{\varepsilon}(t) =
\dot{\mathfrak{X}}(t,s,A_{\varepsilon}) \cdot
\mathfrak{X}^{-1}(t,s,A_{\varepsilon}) - A(t),
\end{equation}
correctness for $\forall t,s \in R$.

Consider the equation
\begin{equation}
\label{kop_eq18} \ddot{y} = \frac{2t^2-1}{(1+t^2)^2} \cdot y.
\end{equation}
For the equation (\ref{kop_eq18}) the equivalence system is
\begin{equation}
\label{kop_eq19} \dot{x} = A(t)\cdot x, \; A(t) = \begin{pmatrix}
0 & 1
\\
\frac{2t^2-1}{(1+t^2)^2} & 0
\end{pmatrix}, \; x \in R^2.
\end{equation}
The Cauchy matrix of system (\ref{kop_eq19}) has a representation:
\begin{equation}
\label{kop_eq20}
\mathfrak{X}(t,s,A)=\frac{1}{3\sqrt{(1+t^2)(1+s^2)}} \cdot
\begin{pmatrix} \frac{f_1(t,s)}{1+s^2} & -3s-s^3+3t+t^3
\\ \frac{f_2(t,s)}{(1+t^2)(1+s^2)} & \frac{3+3ts+ts^3+3t^2+2t^4}{1+t^2}
\end{pmatrix},
\end{equation}
where
$$f_1(t,s)=3+3s^2+2s^4+3ts+t^3s,$$
$$f_2(t,s)=-3t-3ts^2-2ts^4+3s+3st^2+2st^4.$$

Assign for $\forall m \in \{1, \ldots, \overline{t}\}$
($\overline{t} \in N$)
\begin{equation}
\label{kop_eq21} W_m = \mathfrak{X}(m,m-1,A) + \begin{pmatrix} 0 &
0
\\ 0 & \frac{r_{m}}{3\sqrt{(1+m^2)^3(1+(m-1)^2)}}
\end{pmatrix},
\end{equation}
where $r_m > 0$ - some real numbers, that depend on $\delta$ and
$m$.

We have
\begin{equation}
\label{kop_eq22} W_m \cdot [\mathfrak{X}(m,m-1,A)]^{-1} =
\frac{r_{m}}{g_1(m,m-1)} \cdot
\begin{pmatrix} 0 & 0
\\ -f_2(m,m-1) & f_1(m,m-1)
\end{pmatrix},
\end{equation}
\begin{equation}
\label{kop_eq23} \mathfrak{X}(m,m-1,A) \cdot W_m^{-1} =
\frac{r_{m}}{g_2(m,m-1)} \cdot
\begin{pmatrix} 0 & 0
\\ f_2(m,m-1) & f_1(m,m-1)
\end{pmatrix},
\end{equation}
where
$$g_1(t,s)=9 \cdot(1+t^2)^3 \cdot (1+s^2)^2,$$
$$g_2(t,s)=(1+t^2)\cdot (g_3(t,s)+f_1(t,s) \cdot r),$$
$$g_3(t,s)=9(1+2t^2+t^4+4s^2t^2+2s^2t^4+2s^4t^2+s^4t^4+s^2+s^4).$$

If now assign, that the $r_m$ satisfies the inequality
$$
r_m < \frac{\delta}{2\max\{ |f_2(m,m-1)|,|f_1(m,m-1)| \}} \times
$$
$$
\times \max\{ g_1(m,m-1), \; (1+m^2) \cdot (g_3(m,m-1) +
|f_1(m,m-1)| \cdot r_m)\},
$$
then obtain, that the operators $W_m$ for $\forall m \in \{1,
\ldots, \overline{t}\}$ satisfy the inequality (\ref{kop_eq16}).

Then in suppose, for all $m \in \{1, \ldots, \overline{t}\}$
$$
\mathfrak{X}(m,m-1,A_{\varepsilon}) = W_m
$$
and in view of (\ref{kop_eq17}), (\ref{kop_eq21}) and
(\ref{kop_eq20}) we have
$$
B_{\varepsilon}(m)=\dot{W}_m \cdot W_m^{-1} - A(m)=
$$
\begin{equation}
\label{kop_eq24} =\frac{r_m}{g_2(m,m-1)} \cdot
\begin{pmatrix}
f_2(m,m-1) &  -f_1(m,m-1)\cdot (1+m^2) \\
\frac{3f_2(m,m-1) \cdot m}{1+m^2} & 3f_1(m,m-1) \cdot m
\end{pmatrix}.
\end{equation}
In view of (\ref{kop_eq13}), the matrix $B_{\varepsilon}(t)$ must
has a form (\ref{kop_eq15}). But the right side of
(\ref{kop_eq24}) don't satisfies the (\ref{kop_eq15}). This
contradiction proves a lemma.

Consider the family of morphisms $\mathfrak{G}$
\begin{equation}
\label{kop_eq25} (X(m), \chi(m)): \quad (E,p,B) \to (E,p,B),
\end{equation}
$(m \in N)$, of vector bundle $(E,p,B)$, in which connection
\begin{equation}
\label{kop_eq26} B=M_2, \quad E = B \times R^2, \quad p = pr_1,
\end{equation}
\begin{equation}
\label{kop_eq27} X^t(A,x) = (\chi^tA, \mathfrak{X}(t,0,A) \cdot
x),
\end{equation}
\begin{equation}
\label{kop_eq28} \chi^t A(\cdot) = A(t+(\cdot)),
\end{equation}
where $A \in B$, $x \in R^2$, $\mathfrak{X}(\Theta, \tau, A)$ -
Cauchy matrix of the system (\ref{kop_eq12}), that equivalents to
equation (\ref{kop_eq11}).

It is easy to make sure that the family of morphisms
(\ref{kop_eq25}), (\ref{kop_eq27}), (\ref{kop_eq28}) satisfies the
conditions a) - c).

\section{Two-dimensional case}
The aim of this section is to prove the following result.
\begin{lem}
The family of morphisms $\mathfrak{G}$ of vector bundle
(\ref{kop_eq26}) is not saturated.
\end{lem}
{\it Proof}. We argue by contradiction. Suppose, that the family
of morphisms (\ref{kop_eq25}), (\ref{kop_eq27}), (\ref{kop_eq28})
of vector bundle (\ref{kop_eq26}) is saturated.

Then in view of definition, for any point $b \in B$ such that
$\chi^m b \neq b$ for  $\forall m \neq 0$, for $\forall
\overline{\varepsilon}> 0$ for any basis $\{ \xi_1, \xi_2\}$ of
vector space $p^{-1}(b)$ and for all neighbourhoods $U(\xi_i)$ of
points $\xi_i$ $(i=\overline{1,2})$(in space $E$) there exists
$\delta>0$ such that for $\forall \overline{t} \in N$ and for all
nonsingular linear operators $Y_m: \; p^{-1}(\chi^{m-1}b) \to
p^{-1}(\chi^{m}b)$ ($m \in \{ 1, \ldots , \overline{t}\}$), that
satisfy for $\forall m \in \{ 1, \ldots , \overline{t}\}$ the
inequality (\ref{kop_eq4}) there exists the point $b' \in B$ such
that $d_B(b,b')<\overline{\varepsilon}$ and for $\forall m \in \{
0, \ldots , \overline{t}\}$ there exist isomorphisms of fibres (as
a Cartesian spaces) $\psi_m: \; p^{-1}(\chi^{m}b') \to
p^{-1}(\chi^{m}b)$, in which connection, have places conditions i)
- ii).

Let $\overline{\varepsilon} = \varepsilon$ and select $\delta
> 0$ and $\overline{t}$ such that has place the condition of lemma \ref{kop_lm1}.

Fix in the space $B$ any point $b=A$, and then select in the space
$p^{-1}(b)$ any basis $\{\xi_1, \xi_2\}$.

Consider the map $Y_m: \; p^{-1}(\chi^{m-1}b) \to
p^{-1}(\chi^{m}b)$
\begin{equation}
\label{kop_eq32} Y_m(A,\xi)=(\chi^m A,\, W_m \, \xi),
\end{equation}
where $A \in B$, $\xi \in R^2$, $W_m = W(m,m-1)$ ($m \in \{ 1,
\ldots , \overline{t}\}$) - any nonsingular operators from lemma
\ref{kop_lm1}.

It is easy to make sure, that in such case the operator $Y_m$
satisfies the inequality (\ref{kop_eq4}).

Then, in view of supposition, there exists the point $b' \in B$
such that $d_B (b,b') < \overline{\varepsilon}$ and for $\forall m
\in \{ 0, \ldots , \overline{t}\}$ there exist isomorphisms of
fibres $\psi_m: \; p^{-1}(\chi^{m}b') \to p^{-1}(\chi^{m}b)$, in
which connection have places conditions i) - ii).

Let $b'=A_{\varepsilon}$, where
$A_{\varepsilon}=A_{\varepsilon}(t)$ - operator from lemma
\ref{kop_lm1}.

Since $\overline{\varepsilon} = \varepsilon$, then condition $d_B
(b,b') < \overline{\varepsilon}$ has place.

The condition ii) supposes, that diagram (\ref{kop_eq6}) is
commutative, whence, taking into consideration (\ref{kop_eq27}),
(\ref{kop_eq32}), obtain:
\begin{equation}
\label{kop_eq29} \begin{CD}
p^{-1}(\chi^{m-1}b') @>{X[\chi^{m-1}b']}>>  p^{-1}(\chi^{m}b')  \\
@VV{p_{2,\chi^{m-1}b'}}V                            @VV{p_{2,\chi^m b'}}V   \\
R^2 @>\mathfrak{X}(m,m-1,A_{\varepsilon})>> R^2
\end{CD}
\end{equation}
\begin{equation}
\label{kop_eq30} \begin{CD}
p^{-1}(\chi^{m-1}b) @>{Y_m}>>  p^{-1}(\chi^{m}b)  \\
@VV{p_{2,\chi^{m-1}b}}V                            @VV{p_{2,\chi^m b}}V   \\
R^2 @>W_m>> R^2
\end{CD}
\end{equation}
where $p_{2,\overline{b}}$ - the contraction on fibre
$p^{-1}(\overline{b})$ of map $pr_2$ ($pr_2$  - the projection of
product by second efficient).

From the commutative of diagrams (\ref{kop_eq6}),
(\ref{kop_eq29}), (\ref{kop_eq30}) follows commutative of diagram
\begin{equation}
\label{kop_eq31} \begin{CD}
R^2 @>{\mathfrak{X}(m,m-1,A_{\varepsilon})}>>  R^2  \\
@V{p_{2,\chi^{m-1}b} \cdot \psi_{m-1}  \cdot
(p_{2,\chi^{m-1}b'}})^{-1}VV
@VV{p_{2,\chi^m b} \cdot \psi_m \cdot (p_{2,\chi^m b'}})^{-1}V   \\
R^2 @>W_m>> R^2
\end{CD}
\end{equation}
From where follows a claim: for $\forall \, \varepsilon >0$
$\exists \, \delta>0$ such that for $\forall \overline{t} \in N$
there always exist nonsingular linear operators $W_m: \; R^2 \to
R^2$ ($m \in \{1, \ldots, \overline{t} \}$), that satisfy for
$\forall m \in \{1, \ldots, \overline{t} \}$ the inequality
(\ref{kop_eq16}), for which there no exists the continues map
$A_{\varepsilon}(\cdot): \; [0,\overline{t}] \to
\mathrm{Hom}(R^2,R^2)$, that satisfies the conditions 1) - 2).
Since obtained claim is contrary to lemma \ref{kop_lm1}, then
lemma is proved.

\section{n-dimensional case}
Now consider the linear differential equations
\begin{equation}
\label{kop_eq33} y^{(n)}+a_1(t)y^{(n-1)} + \ldots + a_n(t)y=0,
\end{equation}
where $\sup\limits_{t}|a_i(t)|<\infty$, ($i=\overline{1,n}$) and
equivalence to them systems of differential equations
\begin{equation}
\label{kop_eq34} \dot{x}=A(t) \cdot x, \quad x \in R^{n}.
\end{equation}
On the set of systems (\ref{kop_eq34}) defines the metric. The
obtained metric space denotes as $M_n$.

The family of morphisms
$$
(X(m), \chi(m)): (E,p,B) \to (E,p,B),
$$
($m \in N$), of vector bundle $(E,p,B)$, where
\begin{equation}
\label{kop_eq35} B=M_n, \quad E = B \times R^n, \quad p = pr_1,
\end{equation}
$$ X^t(A,x) = (\chi^tA, \mathfrak{X}(t,0,A) \cdot x),
$$
$$ \chi^t A(\cdot) = A(t+(\cdot)),$$
denotes as $\mathfrak{S}$  (here $A \in B$, $x \in R^n$,
$\mathfrak{X}(\Theta, \tau, A)$ - Cauchy matrix of system
(\ref{kop_eq34}), that equivalents to equation (\ref{kop_eq33})).
\begin{thm}
\label{kop_1} If the differential equation (\ref{kop_eq33}) is
given and $\mathfrak{X}(\Theta, \tau,A)$ - a Cauchy matrix of
equivalence system (\ref{kop_eq34}), then for $\forall \varepsilon
> 0$ there exists $\delta > 0$ such that for $\forall \overline{t} \in
N$ there always exist nonsingular linear operators $W_{m}:R^n \to
R^n$ ($m \in \{1,\ldots, \overline{t}\}$) that satisfy for
$\forall m \in \{1, \ldots, \overline{t}\}$ the inequality
(\ref{kop_eq16}), for which there no exists continues map
$$
A_{\varepsilon}(\cdot):[0,\overline{t}] \to \mathrm{Hom}(R^n,R^n)
$$
that satisfies the conditions:

I)$\sup\limits_{t \in [0,\overline{t}]}\|A_{\varepsilon}(t)-A(t)\|
< \varepsilon$;

II) $\mathfrak{X}(m,m-1,A_{\varepsilon})=W_{m}$ for $\forall m \in
\{1, \ldots, \overline{t}\}$,

\noindent where $\mathfrak{X}(m,m-1,A_{\varepsilon})$ - Cauchy
matrix of system (\ref{kop_eq34}) that equivalents to equation
(\ref{kop_eq33}).
\end{thm}

{\it Proof}. Consider the equation
\begin{equation}
\label{kop_eq36} y^{(n)}=a(t)y^{(n-2)}, \quad y \in R.
\end{equation}
The system, that equivalents to (\ref{kop_eq36}), have a
representation:
$$
\dot{x} =A(t) \cdot x, \quad x \in R^n,
$$
where
$$
A(t) =
\begin{pmatrix}
0 & 1 & 0 & \ldots & 0 & 0 \\
0 & 0 & 1 & \ldots & 0 & 0 \\
\ldots & \ldots & \ldots & \ldots & \ldots & \ldots\\
0 & 0 & 0 & \ldots & 1 & 0 \\
0 & 0 & 0 & \ldots & 0 & 1 \\
0 & 0 & 0 & \ldots & a(t) & 0 \\
\end{pmatrix}.
$$
Suppose, that has place the claim, that contraries to claim of
theorem. It means, that for $\forall \overline{t} \in N$ for all
linear operators $W_m$, that satisfy for $\forall m \in \{1,
\ldots, \overline{t}\}$ the inequality (\ref{kop_eq16}), there
always exists a matrix $B_{\varepsilon}(t) = A_{\varepsilon}(t) -
A(t)$ with representation
$$
B_{\varepsilon}(t) =
\begin{pmatrix}
0  & \ldots & 0 & 0 \\
0  & \ldots & 0 & 0 \\
\ldots & \ldots & \ldots & \ldots\\
0 & \ldots & b_{\varepsilon}(t) & 0 \\
\end{pmatrix}, \quad \sup\limits_{t}|b_{\varepsilon}(t)| < \varepsilon,
$$
such that conditions I), II) is satisfied.

On the other hand, if in the equation (\ref{kop_eq36}) we make
change $y^{(n-2)}=z$, then we obtain already considered second
order equation
$$
\ddot{z} = a(t) \cdot z, \quad z \in R^2.
$$
For this equation, according to lemma 1, for $\forall \overline{t}
\in N$ there always exist the linear operators $W_m$, that satisfy
for $\forall m \in \{1, \ldots, \overline{t}\}$ the inequality
(\ref{kop_eq16}), for which there no exists the matrix
$B_{\varepsilon}(t)$, that have representation (\ref{kop_eq15}),
such that the conditions 1), 2) is satisfied. This contradiction
proves theorem.

\begin{thm}
The family of morphisms $\mathfrak{S}$ of vector bundle
(\ref{kop_eq35}) is not saturated.
\end{thm}
The proof of theorem is similarly to proof of lemma 2 and bases on
theorem 1 validity.

\bibliographystyle{plain}
\bibliography{nonexistence}
\end{document}